# Proving ζ(2) through an evolution of the Mengoli's series to the set of rational numbers

URIEL VALENTINIS RAMOS
email: valentinisramos@gmail.com

**ABSTRACT.**

The present paper is an evolution of the Mengoli's series to the set of rational numbers, which eventually will allow developing the summation, by limits, obtaining the value of ζ(2); problem which Mengoli himself was the first to postulate.

More specifically in this paper is postulated and demonstrated the function which resolves every infinite sum of the type $\sum_{n=1}^{\infty}\left(\frac{1}{\left(n+\frac{a}{w}\right)\left(n+\frac{b}{w}\right)\ldots\left(n+\frac{z}{w}\right)}\right)$ for all a, b, …, z, w belonging to the set of integers, where a<b<…<z; $\frac{a}{w}, \frac{b}{w}, \ldots, \frac{z}{w} \neq -1, -2, -3, \ldots$

Finally the following limit $\lim_{w \to \infty}\left\{\sum_{n=1}^{\infty}\left(\frac{1}{\left(n-\frac{1}{w}\right)\left(n+\frac{1}{w}\right)}\right)\right\}$ is developed to demonstrate ζ(2).

**Proposition 1.1.**

All $\sum_{n=1}^{\infty} \left(\frac{1}{(n+a)(n+b)}\right) = \frac{1}{(b-a)} \left[\sum_{n=1}^{b} \left(\frac{1}{n}\right) - \sum_{n=1}^{a} \left(\frac{1}{n}\right)\right]$ where $a, b \in \mathbb{Z}$; $a, b \geq 0$; $a \neq b$.

Proof:

$$\sum_{n=1}^{\infty} \left(\frac{1}{(n+a)(n+b)}\right) = \frac{1}{(b-a)} \left[\sum_{n=1}^{\infty} \left(\frac{n+b-(n+a)}{(n+a)(n+b)}\right)\right] = \frac{1}{(b-a)} \left[\sum_{n=1}^{\infty} \left(\frac{1}{(n+a)}\right) - \sum_{n=1}^{\infty} \left(\frac{1}{(n+b)}\right)\right] =$$

$$= \frac{1}{(b-a)} \left[\left(\sum_{n=1}^{\infty} \left(\frac{1}{n}\right) - \sum_{n=1}^{a} \left(\frac{1}{n}\right)\right) - \left(\sum_{n=1}^{\infty} \left(\frac{1}{n}\right) - \sum_{n=1}^{b} \left(\frac{1}{n}\right)\right)\right] = \frac{1}{(b-a)} \left[\sum_{n=1}^{b} \left(\frac{1}{n}\right) - \sum_{n=1}^{a} \left(\frac{1}{n}\right)\right]$$

**Proposition 1.2.**

All $\sum_{n=1}^{\infty} \left(\frac{1}{\left(n+a+\frac{x}{w}\right)\left(n+b+\frac{y}{w}\right)}\right) =$

$$= \frac{w}{(x+wa-y-wb)} \left[\frac{\pi}{2}\left(\cot\left(\frac{\pi y}{w}\right) - \cot\left(\frac{\pi x}{w}\right)\right) + \sum_{n=1}^{w-1} \left(\ln\left(2\sin\left(\frac{\pi n}{w}\right)\right)\left[\cos\left(\frac{2\pi n x}{w}\right) - \cos\left(\frac{2\pi n y}{w}\right)\right]\right) + \frac{w}{x} - \right.$$

$$\left. \frac{w}{y} + \frac{|a|}{a^{(-1)^{2|a|}}} \sum_{n=\frac{1+\frac{|a|}{a}}{2}}^{|a|-\frac{1-\frac{|a|}{a}}{2}} \left(\frac{w}{x+\frac{|a|}{a}wn}\right) - \frac{|b|}{b^{(-1)^{2|b|}}} \sum_{n=\frac{1+\frac{|b|}{b}}{2}}^{|b|-\frac{1-\frac{|b|}{b}}{2}} \left(\frac{w}{y+\frac{|b|}{b}wn}\right)\right]$$

where $x, y, w \in \mathbb{N}$; $0 < x, y < w$; $a, b \in \mathbb{Z}$; $a + \frac{x}{w} \neq b + \frac{y}{w}$.

Proof:

From Gauss's Theorem for Digamma Function (studied by Jensen [2]) we know that

$$\sum_{n=0}^{\infty} \left(\frac{1}{n+1} - \frac{q}{p+nq}\right) = \frac{-\pi}{2}\cot\left(\frac{\pi p}{q}\right) - \ln(q) + \sum_{n=1}^{q-1} \left(\ln\left(2\sin\left(\frac{\pi n}{q}\right)\right)\cos\left(\frac{2\pi n p}{q}\right)\right)$$

where $p, q \in \mathbb{N}$; $p < q$.

Increasing $n$ from 0 to 1, the following summation is obtained:

$$\sum_{n=1}^{\infty} \left(\frac{1}{n+1} - \frac{q}{p+nq}\right) = \frac{-\pi}{2}\cot\left(\frac{\pi p}{q}\right) - \ln(q) + \sum_{n=1}^{q-1} \left(\ln\left(2\sin\left(\frac{\pi n}{q}\right)\right)\cos\left(\frac{2\pi n p}{q}\right)\right) - 1 + \frac{q}{p}$$

From here two summations are generated:

**Summation 1.2.1**

If $a > 0$

$$\sum_{n=1}^{\infty}\left(\frac{1}{n+1} - \frac{w}{x+w(n+a)}\right) = \sum_{n=1}^{\infty}\left(\frac{1}{n+1}\right) - \sum_{n=1}^{\infty}\left(\frac{w}{x+wn}\right) + \sum_{n=1}^{a}\left(\frac{w}{x+wn}\right)$$

If $a = 0$

$$\sum_{n=1}^{\infty}\left(\frac{1}{n+1} - \frac{w}{x+w(n+a)}\right) = \sum_{n=1}^{\infty}\left(\frac{1}{n+1}\right) - \sum_{n=1}^{\infty}\left(\frac{w}{x+wn}\right)$$

If $a < 0$

$$\sum_{n=1}^{\infty}\left(\frac{1}{n+1} - \frac{w}{x+w(n+a)}\right) = \sum_{n=1}^{\infty}\left(\frac{1}{n+1}\right) - \sum_{n=1}^{\infty}\left(\frac{w}{x+wn}\right) - \sum_{n=0}^{|a|-1}\left(\frac{w}{x-wn}\right)$$

Unifying the three possibilities:

$$\sum_{n=1}^{\infty}\left(\frac{1}{n+1} - \frac{w}{x+w(n+a)}\right) = \sum_{n=1}^{\infty}\left(\frac{1}{n+1} - \frac{w}{x+wn}\right) + \frac{|a|}{a^{(-1)^{2|a|}}}\sum_{n=\frac{1+\frac{|a|}{a}}{2}}^{|a|-\frac{1-\frac{|a|}{a}}{2}}\left(\frac{w}{x+\frac{|a|}{a}wn}\right) =$$

$$= \frac{-\pi}{2}\cot\left(\frac{\pi x}{w}\right) - \ln(w) + \sum_{n=1}^{w-1}\left(\ln\left(2\sin\left(\frac{\pi n}{w}\right)\right)\cos\left(\frac{2\pi n x}{w}\right)\right) - 1 + \frac{w}{x} + \frac{|a|}{a^{(-1)^{2|a|}}}\sum_{n=\frac{1+\frac{|a|}{a}}{2}}^{|a|-\frac{1-\frac{|a|}{a}}{2}}\left(\frac{w}{x+\frac{|a|}{a}wn}\right) \quad (1.2.1)$$

**Summation 1.2.2**

$$\sum_{n=1}^{\infty}\left(\frac{1}{n+1} - \frac{w}{y+w(n+b)}\right) = \sum_{n=1}^{\infty}\left(\frac{1}{n+1} - \frac{w}{y+wn}\right) + \frac{|b|}{b^{(-1)^{2|b|}}}\sum_{n=\frac{1+\frac{|b|}{b}}{2}}^{|b|-\frac{1-\frac{|b|}{b}}{2}}\left(\frac{w}{y+\frac{|b|}{b}wn}\right) =$$

$$= \frac{-\pi}{2}\cot\left(\frac{\pi y}{w}\right) - \ln(w) + \sum_{n=1}^{w-1}\left(\ln\left(2\sin\left(\frac{\pi n}{w}\right)\right)\cos\left(\frac{2\pi n y}{w}\right)\right) - 1 + \frac{w}{y} + \frac{|b|}{b^{(-1)^{2|b|}}}\sum_{n=\frac{1+\frac{|b|}{b}}{2}}^{|b|-\frac{1-\frac{|b|}{b}}{2}}\left(\frac{w}{y+\frac{|b|}{b}wn}\right) \quad (1.2.2)$$

From these expressions, the following operation can be defined **(1.2.1 – 1.2.2):**

$$\sum_{n=1}^{\infty}\left(\frac{1}{n+1}-\frac{w}{x+w(n+a)}\right)-\sum_{n=1}^{\infty}\left(\frac{1}{n+1}-\frac{w}{y+w(n+b)}\right)=\sum_{n=1}^{\infty}\left(\frac{w}{y+w(n+b)}-\frac{w}{x+w(n+a)}\right)=$$

$$=\sum_{n=1}^{\infty}\left(\frac{w(x+wn+wa)-w(y+wn+wb)}{(x+wn+wa)(y+wn+wb)}\right)=\sum_{n=1}^{\infty}\left(\frac{w(x+wa-y-wb)}{(x+wn+wa)(y+wn+wb)}\right)$$

$$=\sum_{n=1}^{\infty}\left(\frac{w(x+wa-y-wb)/w^2}{(x+wn+wa)(y+wn+wb)/w^2}\right)=\sum_{n=1}^{\infty}\left(\frac{(x+wa-y-wb)/w}{\frac{(x+wn+wa)}{w}\frac{(y+wn+wb)}{w}}\right)=$$

$$=\frac{(x+wa-y-wb)}{w}\sum_{n=1}^{\infty}\left(\frac{1}{\frac{(x+wn+wa)}{w}\frac{(y+wn+wb)}{w}}\right)=\frac{(x+wa-y-wb)}{w}\sum_{n=1}^{\infty}\left(\frac{1}{\left(n+a+\frac{x}{w}\right)\left(n+b+\frac{y}{w}\right)}\right)$$

Using this last expression and substituting **1.2.1** and **1.2.2**, the following equation is obtained:

$$\frac{(x+wa-y-wb)}{w}\sum_{n=1}^{\infty}\left(\frac{1}{\left(n+a+\frac{x}{w}\right)\left(n+b+\frac{y}{w}\right)}\right)=$$

$$=\frac{-\pi}{2}\cot\left(\frac{\pi x}{w}\right)-\ln(w)+\sum_{n=1}^{w-1}\left(\ln\left(2\,\sin\left(\frac{\pi n}{w}\right)\right)\cos\left(\frac{2\pi nx}{w}\right)\right)-1+\frac{w}{x}+\frac{|a|}{a^{(-1)^{2|a|}}}\sum_{n=\frac{1+\frac{|a|}{a}}{2}}^{|a|-\frac{1-\frac{|a|}{a}}{2}}\left(\frac{w}{x+\frac{|a|}{a}wn}\right)$$

$$-\left[\frac{-\pi}{2}\cot\left(\frac{\pi y}{w}\right)-\ln(w)+\sum_{n=1}^{w-1}\left(\ln\left(2\,\sin\left(\frac{\pi n}{w}\right)\right)\cos\left(\frac{2\pi ny}{w}\right)\right)-1+\frac{w}{y}+\frac{|b|}{b^{(-1)^{2|b|}}}\sum_{n=\frac{1+\frac{|b|}{b}}{2}}^{|b|-\frac{1-\frac{|b|}{b}}{2}}\left(\frac{w}{y+\frac{|b|}{b}wn}\right)\right]=$$

$$=\frac{w}{(x+wa-y-wb)}\left[\frac{\pi}{2}\left(\cot\left(\frac{\pi y}{w}\right)-\cot\left(\frac{\pi x}{w}\right)\right)+\sum_{n=1}^{w-1}\left(\ln\left(2\sin\left(\frac{\pi n}{w}\right)\right)\left[\cos\left(\frac{2\pi nx}{w}\right)-\cos\left(\frac{2\pi ny}{w}\right)\right]\right)+\frac{w}{x}-\frac{w}{y}+\frac{|a|}{a^{(-1)^{2|a|}}}\sum_{n=\frac{1+\frac{|a|}{a}}{2}}^{|a|-\frac{1-\frac{|a|}{a}}{2}}\left(\frac{w}{x+\frac{|a|}{a}wn}\right)\right.$$

$$\left.-\frac{|b|}{b^{(-1)^{2|b|}}}\sum_{n=\frac{1+\frac{|b|}{b}}{2}}^{|b|-\frac{1-\frac{|b|}{b}}{2}}\left(\frac{w}{y+\frac{|b|}{b}wn}\right)\right]$$

**Proposition 1.3.**

$$\text{All} \sum_{n=1}^{\infty}\left(\frac{1}{\left(n+a+\frac{x}{w}\right)(n+b)}\right) =$$

$$= \frac{w}{(x+wa-wb)}\left[\frac{-\pi}{2}\cot\left(\frac{\pi x}{w}\right) - \ln(w) + \sum_{n=1}^{w-1}\left(\ln\left(2\sin\left(\frac{\pi n}{w}\right)\right)\cos\left(\frac{2\pi nx}{w}\right)\right) - 1 + \frac{w}{x} + \frac{|a|}{a}\sum_{n=\frac{1+\frac{|a|}{a}}{2}}^{|a|-\frac{1-\frac{|a|}{a}}{2}}\left(\frac{w}{x+\frac{|a|}{a}wn}\right)\right.$$

$$\left. + \frac{1}{b+1} - \sum_{n=1}^{b}\left(\frac{1}{n}\right)\right]$$

where $x, w \in \mathbb{N};\ x < w;\ a, b \in \mathbb{Z};\ b \geq 0.$

Proof:

Using **Proposition 1**, the next summation is created:

$$\sum_{n=1}^{\infty}\left(\frac{1}{(n+a)(n+b)}\right) = \frac{1}{(b+1-b)}\left[\sum_{n=1}^{b+1}\left(\frac{1}{n}\right) - \sum_{n=1}^{b}\left(\frac{1}{n}\right)\right] = \frac{1}{b+1} \qquad (1.1)$$

Using summation **1.2.1** (developed in **Proposition 1.2**) and adding the summation **1.1**, the following operation is presented:

$$\sum_{n=1}^{\infty}\left(\frac{1}{n+1} - \frac{w}{x+w(n+a)}\right) + \sum_{n=1}^{\infty}\left(\frac{1}{(n+b)(n+b+1)}\right) = \sum_{n=1}^{\infty}\left(\frac{1}{n+1} - \frac{w}{x+w(n+a)}\right) + \frac{1}{b+1} =$$

$$= \sum_{n=1}^{\infty}\left(\frac{1}{n+1} - \frac{w}{x+w(n+a)}\right) + \left[\sum_{n=1}^{\infty}\left(\frac{1}{n+b}\right) - \sum_{n=1}^{\infty}\left(\frac{1}{n+b+1}\right)\right] = \sum_{n=1}^{\infty}\left(\frac{1}{n+1} - \frac{w}{x+w(n+a)} + \frac{1}{n+b} - \frac{1}{n+b+1}\right) =$$

$$= \sum_{n=1}^{\infty}\left(\frac{1}{n+1} - \frac{1}{n+b+1}\right) + \sum_{n=1}^{\infty}\left(\frac{1}{n+b} - \frac{w}{x+w(n+a)}\right) = \sum_{n=1}^{b}\left(\frac{1}{n+1}\right) + \sum_{n=1}^{\infty}\left(\frac{(x+wn+wa)-w(n+b)}{(x+wn+wa)(n+b)}\right) =$$

$$= \sum_{n=1}^{b}\left(\frac{1}{n+1}\right) + \sum_{n=1}^{\infty}\left(\frac{(x+wa-wb)}{(x+wn+wa)(n+b)}\right) = \sum_{n=1}^{b}\left(\frac{1}{n+1}\right) + \sum_{n=1}^{\infty}\left(\frac{(x+wa-wb)/w}{(x+wn+wa)(n+b)/w}\right) =$$

$$= \sum_{n=1}^{b}\left(\frac{1}{n+1}\right) + \left(\frac{x}{w}+a-b\right)\sum_{n=1}^{\infty}\left(\frac{1}{(n+a+\frac{x}{w})(n+b)}\right)$$

hence,

$$\left(\frac{x}{w} + a - b\right) \sum_{n=1}^{\infty} \left(\frac{1}{(n+a+\frac{x}{w})(n+b)}\right) = \sum_{n=1}^{\infty} \left(\frac{1}{n+1} - \frac{w}{x+w(n+a)}\right) + \sum_{n=1}^{\infty} \left(\frac{1}{(n+b)(n+b+1)}\right) - \sum_{n=1}^{b}\left(\frac{1}{n+1}\right)$$

and substituting both summation **1.1** and **1.2.1**, the next expression is achieved:

$$\left(\frac{x}{w} + a - b\right) \sum_{n=1}^{\infty} \left(\frac{1}{(n+a+\frac{x}{w})(n+b)}\right) =$$

$$= \frac{-\pi}{2}\cot\left(\frac{\pi x}{w}\right) - \ln(w) + \sum_{n=1}^{w-1}\left(\ln\left(2\sin\left(\frac{\pi n}{w}\right)\right)\left[\cos\left(\frac{2\pi n x}{w}\right)\right]\right) - 1 + \frac{w}{x} + \frac{|a|}{a^{(-1)^{2|a|}}} \sum_{n=\frac{1+\frac{|a|}{a}}{2}}^{|a|-\frac{1-\frac{|a|}{a}}{2}} \left(\frac{w}{x+\frac{|a|}{a}wn}\right) + \frac{1}{b+1} - \sum_{n=1}^{b}\left(\frac{1}{n+1}\right)$$

Therefore,

$$\sum_{n=1}^{\infty}\left(\frac{1}{(n+a+x/w)(n+b)}\right) =$$

$$= \frac{w}{(x+wa-wb)}\left[\frac{-\pi}{2}\cot\left(\frac{\pi x}{w}\right) - \ln(w) + \sum_{n=1}^{w-1}\left(\ln\left(2\sin\left(\frac{\pi n}{w}\right)\right)\left[\cos\left(\frac{2\pi n x}{w}\right)\right]\right) + \frac{w}{x} + \frac{|a|}{a^{(-1)^{2|a|}}} \sum_{n=\frac{1+\frac{|a|}{a}}{2}}^{|a|-\frac{1-\frac{|a|}{a}}{2}} \left(\frac{w}{x+\frac{|a|}{a}wn}\right) - 1 + \frac{1}{b+1} - \sum_{n=1}^{b}\left(\frac{1}{n+1}\right)\right] =$$

$$= \frac{w}{(x+wa-wb)}\left[\frac{-\pi}{2}\cot\left(\frac{\pi x}{w}\right) - \ln(w) + \sum_{n=1}^{w-1}\left(\ln\left(2\sin\left(\frac{\pi n}{w}\right)\right)\left[\cos\left(\frac{2\pi n x}{w}\right)\right]\right) + \frac{w}{x} + \frac{|a|}{a^{(-1)^{2|a|}}} \sum_{n=\frac{1+\frac{|a|}{a}}{2}}^{|a|-\frac{1-\frac{|a|}{a}}{2}} \left(\frac{w}{x+\frac{|a|}{a}wn}\right) - \sum_{n=1}^{b}\left(\frac{1}{n}\right)\right]$$

## Proposition 1.4.

All $\sum_{n=1}^{\infty}\left(\frac{1}{\left(n+a+\frac{x}{w}\right)\left(n+b+\frac{y}{w}\right)}\right) =$

$= \frac{w}{(x+wa-y-wb)}\left[\frac{\pi}{2}\left(\frac{\cos\left(\frac{\pi y}{w}\right)}{\left[\sin\left(\frac{\pi y}{w}\right)\right]^{(-1)^{2y}}} - \frac{\cos\left(\frac{\pi x}{w}\right)}{\left[\sin\left(\frac{\pi x}{w}\right)\right]^{(-1)^{2x}}}\right) - \frac{x}{x^{(-1)^{2x}}}\ln(w) + \frac{y}{y^{(-1)^{2y}}}\ln(w)\right.$

$+ \sum_{n=1}^{w-1}\left(\ln\left(2\sin\left(\frac{\pi n}{w}\right)\right)\left[\frac{x}{x^{(-1)^{2x}}}\cos\left(\frac{2\pi nx}{w}\right) - \frac{y}{y^{(-1)^{2y}}}\cos\left(\frac{2\pi ny}{w}\right)\right]\right) + \frac{w}{x^{(-1)^{2x}}} - \frac{w}{y^{(-1)^{2y}}}$

$+ \frac{|a|}{a^{(-1)^{2|a|}}}\sum_{n=\frac{1+\frac{|a|}{a}}{2}}^{|a|-\frac{1-\frac{|a|}{a}}{2}}\left(\frac{w}{x+\frac{|a|}{a}wn}\right) - \frac{|b|}{b^{(-1)^{2|b|}}}\sum_{n=\frac{1+\frac{|b|}{b}}{2}}^{|b|-\frac{1-\frac{|b|}{b}}{2}}\left(\frac{w}{y+\frac{|b|}{b}wn}\right)\right] =$

where $a, b, x, y, w \in \mathbb{Z}$; $0 \leq x, y < w$; $a + \frac{x}{w} \neq b + \frac{y}{w}$; $\left(a + \frac{x}{w}\right), \left(b + \frac{y}{w}\right) \neq -1, -2, -3, \ldots$

Proof:

**Proposition 1.4** is obtained by consolidation of Propositions 1.1, 1.2 and 1.3.

- From Proposition 1.1: All $\sum_{n=1}^{\infty}\left(\frac{1}{(n+a)(n+b)}\right) = \frac{1}{(b-a)}\left[\sum_{n=1}^{b}\left(\frac{1}{n}\right) - \sum_{n=1}^{a}\left(\frac{1}{n}\right)\right]$
  where $a, b \in \mathbb{Z}$; $a, b \geq 0$; $a \neq b$.
  And $\sum_{n=1}^{\infty}\left(\frac{1}{(n+a)(n+b)}\right) = \sum_{n=1}^{\infty}\left(\frac{1}{\left(n+a+\frac{x}{w}\right)\left(n+b+\frac{y}{w}\right)}\right) \Rightarrow x, y = 0, w \neq 0$

Substituting $x, y = 0$; $w \neq 0$; $a, b \geq 0$ at **Proposition 1.4**:

$$= \frac{w}{(0 + wa - 0 - wb)} \left[ \frac{\pi}{2} \left( \frac{\cos\left(\frac{\pi 0}{w}\right)}{\left[\sin\left(\frac{\pi 0}{w}\right)\right]^{(-1)^{2^0}}} - \frac{\cos\left(\frac{\pi 0}{w}\right)}{\left[\sin\left(\frac{\pi 0}{w}\right)\right]^{(-1)^{2^0}}} \right) - \frac{0}{0^{(-1)^{2|0|}}} \ln(w) + \frac{0}{0^{(-1)^{2|0|}}} \ln(w) \right.$$

$$+ \sum_{n=1}^{w-1} \left( \ln\left(2 \sin\left(\frac{\pi n}{w}\right)\right) \left[ \frac{0}{0^{(-1)^{2|0|}}} \cos\left(\frac{2\pi n 0}{w}\right) - \frac{0}{0^{(-1)^{2|0|}}} \cos\left(\frac{2\pi n 0}{w}\right) \right] \right) + \frac{w}{0^{(-1)^{2|0|}}} - \frac{w}{0^{(-1)^{2|0|}}}$$

$$\left. + \frac{|a|}{+a^{(-1)^{2|a|}}} \sum_{n=\frac{1+\frac{|a|}{a}}{2}}^{|+a|-\frac{1-\frac{|a|}{a}}{2}} \left( \frac{\cancel{w}}{0 + \frac{|a|}{+a} \cancel{w} n} \right) - \frac{|b|}{b^{(-1)^{2|b|}}} \sum_{n=\frac{1+\frac{|b|}{b}}{2}}^{|b|-\frac{1-\frac{|b|}{b}}{2}} \left( \frac{\cancel{w}}{0 + \frac{|b|}{b} \cancel{w} n} \right) \right] =$$

$$= \frac{\cancel{w}}{\cancel{w}(a-b)} \left[ \sum_{n=1}^{a} \left(\frac{1}{n}\right) - \sum_{n=1}^{b} \left(\frac{1}{n}\right) \right] = \frac{1}{(a-b)} \left[ \sum_{n=1}^{a} \left(\frac{1}{n}\right) - \sum_{n=1}^{b} \left(\frac{1}{n}\right) \right]$$

- From Proposition 1.2: *All* $\displaystyle\sum_{n=1}^{\infty} \left( \frac{1}{\left(n + a + \frac{x}{w}\right)\left(n + b + \frac{y}{w}\right)} \right)$

$$= \frac{w}{(x + wa - y - wb)} \left[ \frac{\pi}{2} \left( \cot\left(\frac{\pi y}{w}\right) - \cot\left(\frac{\pi x}{w}\right) \right) + \sum_{n=1}^{w-1} \left( \ln\left(2 \sin\left(\frac{\pi n}{w}\right)\right) \left[ \cos\left(\frac{2\pi n x}{w}\right) - \cos\left(\frac{2\pi n y}{w}\right) \right] \right) + \frac{w}{x} - \frac{w}{y} \right.$$

$$\left. + \frac{|a|}{a^{(-1)^{2|a|}}} \sum_{n=\frac{1+\frac{|a|}{a}}{2}}^{|a|-\frac{1-\frac{|a|}{a}}{2}} \left( \frac{w}{x + \frac{|a|}{a} w n} \right) - \frac{|b|}{b^{(-1)^{2|b|}}} \sum_{n=\frac{1+\frac{|b|}{b}}{2}}^{|b|-\frac{1-\frac{|b|}{b}}{2}} \left( \frac{w}{y + \frac{|b|}{b} w n} \right) \right]$$

*where* $x, y, w \in \mathbb{N}$; $0 < x, y < w$; $a, b \in \mathbb{Z}$; $a + \frac{x}{w} \neq b + \frac{y}{w}$.

Substituting $0 < x, y < w$ at **Proposition 1.4**:

$$= \frac{w}{(x+wa-y-wb)}\left[\frac{\pi}{2}\left(\frac{\cos\left(\frac{\pi y}{w}\right)}{\sin\left(\frac{\pi y}{w}\right)} - \frac{\cos\left(\frac{\pi x}{w}\right)}{\sin\left(\frac{\pi y}{w}\right)}\right) - \frac{\cancel{x}}{\cancel{x}}\ln(w) + \frac{\cancel{y}}{\cancel{y}}\ln(w) + \sum_{n=1}^{w-1}\left(\ln\left(2\sin\left(\frac{\pi n}{w}\right)\right)\left[\frac{\cancel{x}}{\cancel{x}}\cos\left(\frac{2\pi nx}{w}\right) - \frac{\cancel{y}}{\cancel{y}}\cos\left(\frac{2\pi ny}{w}\right)\right]\right)\right.$$

$$\left. + \frac{w}{x} - \frac{w}{y} + \frac{|a|}{a^{(-1)^{2|a|}}}\sum_{n=\frac{1+\frac{|a|}{a}}{2}}^{|a|-\frac{1-\frac{|a|}{a}}{2}}\left(\frac{w}{x+\frac{|a|}{a}wn}\right) - \frac{|b|}{b^{(-1)^{2|b|}}}\sum_{n=\frac{1+\frac{|b|}{b}}{2}}^{|b|-\frac{1-\frac{|b|}{b}}{2}}\left(\frac{w}{y+\frac{|b|}{b}wn}\right)\right] =$$

$$= \frac{w}{(x+wa-y-wb)}\left[\frac{\pi}{2}\left(\cot\left(\frac{\pi y}{w}\right) - \cot\left(\frac{\pi x}{w}\right)\right) + \sum_{n=1}^{w-1}\left(\ln\left(2\sin\left(\frac{\pi n}{w}\right)\right)\left[\cos\left(\frac{2\pi nx}{w}\right) - \cos\left(\frac{2\pi ny}{w}\right)\right]\right) + \frac{w}{x} - \frac{w}{y}\right.$$

$$\left. + \frac{|a|}{a^{(-1)^{2|a|}}}\sum_{n=\frac{1+\frac{|a|}{a}}{2}}^{|a|-\frac{1-\frac{|a|}{a}}{2}}\left(\frac{w}{x+\frac{|a|}{a}wn}\right) - \frac{|b|}{b^{(-1)^{2|b|}}}\sum_{n=\frac{1+\frac{|b|}{b}}{2}}^{|b|-\frac{1-\frac{|b|}{b}}{2}}\left(\frac{w}{y+\frac{|b|}{b}wn}\right)\right]$$

- From Proposition 1.3: *All* $\displaystyle\sum_{n=1}^{\infty}\left(\frac{1}{\left(n+a+\frac{x}{w}\right)(n+b)}\right) =$

$$= \frac{w}{(x+wa-wb)}\left[\frac{-\pi}{2}\cot\left(\frac{\pi x}{w}\right) - \ln(w) + \sum_{n=1}^{w-1}\left(\ln\left(2\sin\left(\frac{\pi n}{w}\right)\right)\cos\left(\frac{2\pi nx}{w}\right)\right) - 1 + \frac{w}{x}\right.$$

$$\left. + \frac{|a|}{a}\sum_{n=\frac{1+\frac{|a|}{a}}{2}}^{|a|-\frac{1-\frac{|a|}{a}}{2}}\left(\frac{w}{x+\frac{|a|}{a}wn}\right) + \frac{1}{b+1} - \sum_{n=1}^{b}\left(\frac{1}{n}\right)\right] \quad \text{where} \quad x, w \in \mathbb{N}; \; x < w; \; a, b \in \mathbb{Z}; \; b \geq 0.$$

Substituting $y = 0;\ 0 < x < w;\ b \geq 0$ at **Proposition 1.4**:

$$= \frac{w}{(x + wa - 0 - wb)} \left[ \frac{\pi}{2}\left(0 - \frac{\cos\left(\frac{\pi x}{w}\right)}{\sin\left(\frac{\pi x}{w}\right)}\right) - \frac{\cancel{*}}{\cancel{*}}\ln(w) + 0\ \ln(w) \right.$$

$$+ \sum_{n=1}^{w-1} \left( \ln\left(2\ \sin\left(\frac{\pi n}{w}\right)\right) \left[\frac{\cancel{*}}{\cancel{*}} \cos\left(\frac{2\pi n x}{w}\right) - 0 \cos\left(\frac{2\pi n y}{w}\right)\right]\ \right) + \frac{w}{x} - 0\ w$$

$$\left. + \frac{|a|}{a^{(-1)^{2|a|}}} \sum_{n=\frac{1+\frac{|a|}{a}}{2}}^{|a|-\frac{1-\frac{|a|}{a}}{2}} \left( \frac{w}{x + \frac{|a|}{a} wn} \right) - \sum_{n=1}^{|b|} \left( \frac{\cancel{w}}{0 + \cancel{w}n} \right) \right] =$$

$$= \frac{w}{(x + wa - wb)} \left[ \frac{-\pi}{2}\cot\left(\frac{\pi x}{w}\right) - \ln(w) + \sum_{n=1}^{w-1} \left( \ln\left(2\ \sin\left(\frac{\pi n}{w}\right)\right) \left[\cos\left(\frac{2\pi n x}{w}\right)\right]\ \right) + \frac{w}{x} + \frac{|a|}{a^{(-1)^{2|a|}}} \sum_{n=\frac{1+\frac{|a|}{a}}{2}}^{|a|-\frac{1-\frac{|a|}{a}}{2}} \left( \frac{w}{x + \frac{|a|}{a} wn} \right) \right.$$

$$\left. - \sum_{n=1}^{b} \left(\frac{1}{n}\right) \right]$$

An example of the proposition 1.4 can be found in the development of **ζ(2), which is proved by limits**, as follows:

$$\sum_{n=1}^{\infty} \left( \frac{1}{\left(n + (-1) + \frac{w-1}{w}\right)\left(n + 0 + \frac{1}{w}\right)} \right) =$$

$$= \frac{w}{(w - 1 - w - 1)} \left[ \frac{\pi}{2}\left(\cot\left(\frac{\pi}{w}\right) - \cot\left(\frac{\pi(w-1)}{w}\right)\right) + \sum_{n=1}^{w-1} \left( \ln\left(2\sin\left(\frac{\pi n}{w}\right)\right) \left[\cos\left(\frac{2\pi n(w-1)}{w}\right) - \cos\left(\frac{2\pi n}{w}\right)\right]\ \right) + \frac{w}{(w-1)} - \frac{w}{1} - \frac{w}{(w-1)} \right] =$$

$$= \frac{w}{(-2)} \left[ \frac{\pi}{2}\left(\cot\left(\frac{\pi}{w}\right) - \cot\left(\frac{\pi(w-1)}{w}\right)\right) + \sum_{n=1}^{w-1} \left( \ln\left(2\sin\left(\frac{\pi n}{w}\right)\right) \left[\cos\left(\frac{2\pi n(w-1)}{w}\right) - \cos\left(\frac{2\pi n}{w}\right)\right]\ \right) - \frac{w}{1} \right] =$$

$$= \frac{w}{2} \left[ -\frac{\pi}{2}\left(2 \cot\left(\frac{\pi}{w}\right)\right) - \sum_{n=1}^{w-1} \left( \ln\left(2\sin\left(\frac{\pi n}{w}\right)\right)[0]\ \right) + w \right] = w\left[\frac{w}{2} - \frac{\pi}{2}\cot\left(\frac{\pi}{w}\right)\right]$$

Applying limits:

$$\lim_{w \to \infty} \left\{ \sum_{n=1}^{\infty} \left( \frac{1}{\left(n + (-1) + \frac{w-1}{w}\right)\left(n + \frac{1}{w}\right)} \right) \right\} = \sum_{n=1}^{\infty} \left(\frac{1}{n^2}\right) = \lim_{w \to \infty} \left\{ w \left[\frac{w}{2} - \frac{\pi}{2} \cot\left(\frac{\pi}{w}\right)\right] \right\}$$

Due to $\lim_{w \to \infty} \left\{ w \left[\frac{w}{2} - \frac{\pi}{2} \cot\left(\frac{\pi}{w}\right)\right] \right\} = c$, $f(w) = \left[\frac{w^2}{2} - \frac{\pi w}{2} \cot\left(\frac{\pi}{w}\right)\right]$ being continuous and $\lim_{w \to 0}\{1/w\} = \infty$ these imply $\lim_{w \to \infty} \left\{ w \left[\frac{w}{2} - \frac{\pi}{2} \cot\left(\frac{\pi}{w}\right)\right] \right\} = \lim_{w \to 0} \left\{ 1/w \left[\frac{1/w}{2} - \frac{\pi}{2} \cot\left(\frac{\pi}{1/w}\right)\right] \right\}$

$$\lim_{w \to 0} \left\{ 1/w \left[\frac{1/w}{2} - \frac{\pi}{2} \cot\left(\frac{\pi}{1/w}\right)\right] \right\} = \lim_{w \to 0} \left\{ 1/w \left[\frac{1/w}{2} - \frac{\pi}{2} \cot\left(\frac{\pi}{1/w}\right)\right] \right\} = \lim_{w \to 0} \left\{ \frac{1}{2w^2} \right\} = \lim_{w \to 0} \left\{ \frac{\tan(\pi w) - \pi w}{2w^2 \tan(\pi w)} - \frac{\pi}{2w} \cot(\pi w) \right\} = \frac{0}{0}$$

L'Hopital is applied

$$\lim_{w \to 0} \left\{ \frac{\frac{\pi[\cos^2(\pi w) + \sin^2(\pi w)]}{\cos^2(\pi w)} - \pi}{2\left[2w \tan(\pi w) + \frac{w^2 \pi[\cos^2(\pi w) + \sin^2(\pi w)]}{\cos^2(\pi w)}\right]} \right\} = \pi/2 * \lim_{w \to 0} \left\{ \frac{\frac{1}{\cos^2(\pi w)} - 1}{2w \tan(\pi w) + \frac{w^2 \pi}{\cos^2(\pi w)}} \right\} =$$

$$= \pi/2 * \lim_{w \to 0} \left\{ \frac{\cos^2(\pi w)}{\cos^2(\pi w)} \frac{\frac{1}{\cos^2(\pi w)} - 1}{2w \tan(\pi w) + \frac{w^2 \pi}{\cos^2(\pi w)}} \right\}$$

$$= \pi/2 * \lim_{w \to 0} \left\{ \frac{1 - \cos^2(\pi w)}{2w \sin(\pi w)\cos(\pi w) + w^2 \pi} \right\} = \pi/2 * \lim_{w \to 0} \left\{ \frac{1 - \cos^2(\pi w)}{w \sin(2\pi w) + w^2 \pi} \right\} = \frac{0}{0}$$

L'Hopital is applied

$$\pi/2 * \lim_{w \to 0} \left\{ \frac{-[-2\pi \sin(\pi w)\cos(\pi w)]}{\sin(2\pi w) + 2\pi w \cos(2\pi w) + 2w\pi} \right\} = \pi/2 * \lim_{w \to 0} \left\{ \frac{\pi \sin(2\pi w)}{\sin(2\pi w) + 2\pi w \cos(2\pi w) + 2w\pi} \right\}$$

$$= \pi^2/2 * \lim_{w \to 0} \left\{ \frac{\sin(2\pi w)}{\sin(2\pi w) + 2\pi w \cos(2\pi w) + 2w\pi} \right\} = \frac{0}{0}$$

L'Hopital is applied

$$\pi^2/2 * \lim_{w \to 0} \left\{ \frac{2\pi \cos(2\pi w)}{2\pi \cos(2\pi w) + 2\pi \cos(2\pi w) - 4\pi^2 w \sin(2\pi w) + 2\pi} \right\} =$$

$$= \pi^2/2 * \lim_{w \to 0} \left\{ \frac{\cos(2\pi w)}{\cos(2\pi w) + \cos(2\pi w) - 2\pi w \sin(2\pi w) + 1} \right\} = \pi^2/2 * \left[\frac{1}{1 + 1 - 0 + 1}\right] = \frac{\pi^2}{6}$$

## Proposition 2.

Every infinite sum with z divisors: $\left(n+\frac{a}{w}\right)\left(n+\frac{b}{w}\right)\ldots\left(n+\frac{z}{w}\right)$ can be solved based on a finite set of infinite sums with 2 divisors. This means, that every

$$\sum_{n=1}^{\infty}\left(\frac{1}{\left(n+\frac{a}{w}\right)\left(n+\frac{b}{w}\right)\ldots\left(n+\frac{z}{w}\right)}\right)$$

Can be solved as follows

$$\frac{\sum_{n=1}^{\infty}\left(\frac{1}{\left(n+\frac{a}{w}\right)\ldots\left(n+\frac{z}{w}\right)}\right)-\sum_{n=1}^{\infty}\left(\frac{1}{\left(n+\frac{b}{w}\right)\ldots\left(n+\frac{z}{w}\right)}\right)}{\frac{b}{w}-\frac{a}{w}}$$

Proof:

$$\frac{1}{\left(n+\frac{a}{w}\right)\left(n+\frac{b}{w}\right)\ldots\left(n+\frac{z}{w}\right)}=\frac{1}{\frac{b}{w}-\frac{a}{w}}\left[\frac{1}{\left(n+\frac{a}{w}\right)\ldots\left(n+\frac{z}{w}\right)}-\frac{1}{\left(n+\frac{b}{w}\right)\ldots\left(n+\frac{z}{w}\right)}\right]=$$

$$=\frac{1}{\frac{b}{w}-\frac{a}{w}}\left[\frac{\left(n+\frac{b}{w}\right)-\left(n+\frac{a}{w}\right)}{\left(n+\frac{a}{w}\right)\left(n+\frac{b}{w}\right)\ldots\left(n+\frac{z}{w}\right)}\right]=\frac{1}{\frac{b}{w}-\frac{a}{w}}\left[\frac{\frac{b}{w}-\frac{a}{w}}{\left(n+\frac{a}{w}\right)\left(n+\frac{b}{w}\right)\ldots\left(n+\frac{z}{w}\right)}\right]=$$

$$=\frac{1}{\left(n+\frac{a}{w}\right)\left(n+\frac{b}{w}\right)\ldots\left(n+\frac{z}{w}\right)}$$

An example of the proposition 2 can be found in the development of **ζ(4)**, as follows:

$$\sum_{n=1}^{\infty}\left(\frac{1}{\left(n+\frac{-2}{w}\right)\left(n+\frac{-1}{w}\right)\left(n+\frac{1}{w}\right)\left(n+\frac{2}{w}\right)}\right)=\sum_{n=1}^{\infty}\left(\frac{1}{\left(n+(-1)+\frac{w-2}{w}\right)\left(n+(-1)+\frac{w-1}{w}\right)\left(n+\frac{1}{w}\right)\left(n+\frac{2}{w}\right)}\right)=$$

$$=\frac{\sum_{n=1}^{\infty}\left(\frac{1}{\left(n+(-1)+\frac{w-2}{w}\right)\left(n+\frac{1}{w}\right)\left(n+\frac{2}{w}\right)}\right)-\sum_{n=1}^{\infty}\left(\frac{1}{\left(n+(-1)+\frac{w-1}{w}\right)\left(n+\frac{1}{w}\right)\left(n+\frac{2}{w}\right)}\right)}{\left((-1)+\frac{w-1}{w}\right)-\left((-1)+\frac{w-2}{w}\right)}=$$

$$=\frac{\sum_{n=1}^{\infty}\left(\frac{1}{\left(n+(-1)+\frac{w-2}{w}\right)\left(n+\frac{1}{w}\right)\left(n+\frac{2}{w}\right)}\right)-\sum_{n=1}^{\infty}\left(\frac{1}{\left(n+(-1)+\frac{w-1}{w}\right)\left(n+\frac{1}{w}\right)\left(n+\frac{2}{w}\right)}\right)}{\frac{w-1-w+2}{w}}=$$

$$= w\left[\sum_{n=1}^{\infty}\left(\frac{1}{\left(n+(-1)+\frac{w-2}{w}\right)\left(n+\frac{1}{w}\right)\left(n+\frac{2}{w}\right)}\right) - \sum_{n=1}^{\infty}\left(\frac{1}{\left(n+(-1)+\frac{w-1}{w}\right)\left(n+\frac{1}{w}\right)\left(n+\frac{2}{w}\right)}\right)\right] =$$

$$= w\left[\frac{\sum_{n=1}^{\infty}\left(\frac{1}{\left(n+(-1)+\frac{w-2}{w}\right)\left(n+\frac{2}{w}\right)}\right) - \sum_{n=1}^{\infty}\left(\frac{1}{\left(n+\frac{1}{w}\right)\left(n+\frac{2}{w}\right)}\right)}{\left(\frac{1}{w}\right) - \left((-1)+\frac{w-2}{w}\right)} - \frac{\sum_{n=1}^{\infty}\left(\frac{1}{\left(n+(-1)+\frac{w-1}{w}\right)\left(n+\frac{2}{w}\right)}\right) - \sum_{n=1}^{\infty}\left(\frac{1}{\left(n+\frac{1}{w}\right)\left(n+\frac{2}{w}\right)}\right)}{\left(\frac{1}{w}\right) - \left((-1)+\frac{w-1}{w}\right)}\right] =$$

$$= w\left[\frac{\sum_{n=1}^{\infty}\left(\frac{1}{\left(n+(-1)+\frac{w-2}{w}\right)\left(n+\frac{2}{w}\right)}\right) - \sum_{n=1}^{\infty}\left(\frac{1}{\left(n+\frac{1}{w}\right)\left(n+\frac{2}{w}\right)}\right)}{1+\left(\frac{1-w+2}{w}\right)} - \frac{\sum_{n=1}^{\infty}\left(\frac{1}{\left(n+(-1)+\frac{w-1}{w}\right)\left(n+\frac{2}{w}\right)}\right) - \sum_{n=1}^{\infty}\left(\frac{1}{\left(n+\frac{1}{w}\right)\left(n+\frac{2}{w}\right)}\right)}{1+\frac{1-w+1}{w}}\right] =$$

$$= w\left[\frac{\sum_{n=1}^{\infty}\left(\frac{1}{\left(n+(-1)+\frac{w-2}{w}\right)\left(n+\frac{2}{w}\right)}\right) - \sum_{n=1}^{\infty}\left(\frac{1}{\left(n+\frac{1}{w}\right)\left(n+\frac{2}{w}\right)}\right)}{\frac{3}{w}} - \frac{\sum_{n=1}^{\infty}\left(\frac{1}{\left(n+(-1)+\frac{w-1}{w}\right)\left(n+\frac{2}{w}\right)}\right) - \sum_{n=1}^{\infty}\left(\frac{1}{\left(n+\frac{1}{w}\right)\left(n+\frac{2}{w}\right)}\right)}{\frac{2}{w}}\right] =$$

$$= \frac{w^2}{6}\left[2\sum_{n=1}^{\infty}\left(\frac{1}{\left(n+(-1)+\frac{w-2}{w}\right)\left(n+\frac{2}{w}\right)}\right) - 2\sum_{n=1}^{\infty}\left(\frac{1}{\left(n+\frac{1}{w}\right)\left(n+\frac{2}{w}\right)}\right) - 3\sum_{n=1}^{\infty}\left(\frac{1}{\left(n+(-1)+\frac{w-1}{w}\right)\left(n+\frac{2}{w}\right)}\right) + 3\sum_{n=1}^{\infty}\left(\frac{1}{\left(n+\frac{1}{w}\right)\left(n+\frac{2}{w}\right)}\right)\right] =$$

$$= \frac{w^2}{6}\left[2\sum_{n=1}^{\infty}\left(\frac{1}{\left(n+(-1)+\frac{w-2}{w}\right)\left(n+\frac{2}{w}\right)}\right) + \sum_{n=1}^{\infty}\left(\frac{1}{\left(n+\frac{1}{w}\right)\left(n+\frac{2}{w}\right)}\right) - 3\sum_{n=1}^{\infty}\left(\frac{1}{\left(n+(-1)+\frac{w-1}{w}\right)\left(n+\frac{2}{w}\right)}\right)\right] =$$

$$= \frac{w^2}{3}\left[\sum_{n=1}^{\infty}\left(\frac{1}{\left(n+(-1)+\frac{w-2}{w}\right)\left(n+\frac{2}{w}\right)}\right) - \sum_{n=1}^{\infty}\left(\frac{1}{\left(n+(-1)+\frac{w-1}{w}\right)\left(n+\frac{1}{w}\right)}\right)\right]$$

Proposition 1 is applied

$$\frac{w^2}{3}\left[\sum_{n=1}^{\infty}\left(\frac{1}{\left(n+(-1)+\frac{w-2}{w}\right)\left(n+\frac{2}{w}\right)}\right) - \sum_{n=1}^{\infty}\left(\frac{1}{\left(n+(-1)+\frac{w-1}{w}\right)\left(n+\frac{1}{w}\right)}\right)\right] =$$

$$= \frac{w^2}{6}\left[2\left(\frac{w}{((w-2)-w-2)}\left[\frac{\pi}{2}\left(\cot\left(\frac{2\pi}{w}\right) - \cot\left(\frac{\pi(w-2)}{w}\right)\right) + \sum_{n=1}^{w-1}\left(\ln\left(2\sin\left(\frac{\pi n}{w}\right)\right)\left[\cos\left(\frac{2\pi n(w-2)}{w}\right) - \cos\left(\frac{2\pi n 2}{w}\right)\right]\right) + \frac{w}{(w-2)} - \frac{w}{2} - \frac{w}{(w-2)}\right]\right)$$

$$+ \left(\frac{w}{(1-2)}\left[\frac{\pi}{2}\left(\cot\left(\frac{2\pi}{w}\right) - \cot\left(\frac{\pi}{w}\right)\right) + \sum_{n=1}^{w-1}\left(\ln\left(2\sin\left(\frac{\pi n}{w}\right)\right)\left[\cos\left(\frac{2\pi n}{w}\right) - \cos\left(\frac{2\pi n 2}{w}\right)\right]\right) + \frac{w}{1} - \frac{w}{2}\right]\right)$$

$$- 3\left(\frac{w}{((w-1)-w-2)}\left[\frac{\pi}{2}\left(\cot\left(\frac{2\pi}{w}\right) - \cot\left(\frac{\pi(w-1)}{w}\right)\right) + \sum_{n=1}^{w-1}\left(\ln\left(2\sin\left(\frac{\pi n}{w}\right)\right)\left[\cos\left(\frac{2\pi n(w-1)}{w}\right) - \cos\left(\frac{2\pi n 2}{w}\right)\right]\right) + \frac{w}{(w-1)} - \frac{w}{2}\right.\right.$$

$$\left.\left.- \frac{w}{(w-1)}\right]\right)\right] =$$

$$= \frac{w^2}{6}\left[2\left(\frac{w}{(-4)}\left[\frac{\pi}{2}\left(\cot\left(\frac{2\pi}{w}\right) - \cot\left(\frac{\pi(w-2)}{w}\right)\right) + \sum_{n=1}^{w-1}\left(\ln\left(2\sin\left(\frac{\pi n}{w}\right)\right)[0]\right) - \frac{w}{2}\right]\right)\right.$$

$$+ \left(\frac{w}{(-1)}\left[\frac{\pi}{2}\left(\cot\left(\frac{2\pi}{w}\right) - \cot\left(\frac{\pi}{w}\right)\right) + \sum_{n=1}^{w-1}\left(\ln\left(2\sin\left(\frac{\pi n}{w}\right)\right)\left[\cos\left(\frac{2\pi n}{w}\right) - \cos\left(\frac{2\pi n 2}{w}\right)\right]\right) + \frac{w}{2}\right]\right)$$

$$\left.- 3\left(\frac{w}{(-3)}\left[\frac{\pi}{2}\left(\cot\left(\frac{2\pi}{w}\right) - \cot\left(\frac{\pi(w-1)}{w}\right)\right) + \sum_{n=1}^{w-1}\left(\ln\left(2\sin\left(\frac{\pi n}{w}\right)\right)\left[\cos\left(\frac{2\pi n(w-1)}{w}\right) - \cos\left(\frac{2\pi n 2}{w}\right)\right]\right) - \frac{w}{2}\right]\right)\right] =$$

$$= \frac{w^3}{6}\left[\frac{\pi}{4}\left(\cot\left(\frac{\pi(w-2)}{w}\right) - \cot\left(\frac{2\pi}{w}\right)\right) + \frac{w}{4} - \frac{\pi}{2}\left(\cot\left(\frac{2\pi}{w}\right) - \cot\left(\frac{\pi}{w}\right)\right) - \sum_{n=1}^{w-1}\left(\ln\left(2\sin\left(\frac{\pi n}{w}\right)\right)\left[\cos\left(\frac{2\pi n}{w}\right) - \cos\left(\frac{2\pi n 2}{w}\right)\right]\right) - \frac{w}{2}\right.$$

$$\left.+ \frac{\pi}{2}\left(\cot\left(\frac{2\pi}{w}\right) - \cot\left(\frac{\pi(w-1)}{w}\right)\right) + \sum_{n=1}^{w-1}\left(\ln\left(2\sin\left(\frac{\pi n}{w}\right)\right)\left[\cos\left(\frac{2\pi n(w-1)}{w}\right) - \cos\left(\frac{2\pi n 2}{w}\right)\right]\right) - \frac{w}{2}\right] =$$

$$= \frac{w^3}{6}\left[\frac{\pi}{4}\left(\cot\left(\frac{\pi(w-2)}{w}\right) - \cot\left(\frac{2\pi}{w}\right) + 2\cot\left(\frac{\pi}{w}\right) - 2\cot\left(\frac{\pi(w-1)}{w}\right)\right) + \frac{w}{4} - w\right.$$

$$\left.+ \sum_{n=1}^{w-1}\left(\ln\left(2\sin\left(\frac{\pi n}{w}\right)\right)\left[\cos\left(\frac{2\pi n(w-1)}{w}\right) - \cos\left(\frac{2\pi n 2}{w}\right) - \cos\left(\frac{2\pi n}{w}\right) + \cos\left(\frac{2\pi n 2}{w}\right)\right]\right)\right] =$$

$$= \frac{w^3}{6}\left[\frac{\pi}{4}\left(4\cot\left(\frac{\pi}{w}\right) - 2\cot\left(\frac{2\pi}{w}\right)\right) - \frac{3w}{4} + \sum_{n=1}^{w-1}\left(\ln\left(2\sin\left(\frac{\pi n}{w}\right)\right)[0]\right)\right] = \frac{w^3}{24}\left[\pi\left(4\cot\left(\frac{\pi}{w}\right) - 2\cot\left(\frac{2\pi}{w}\right)\right) - 3w\right]$$

Finally, applying limits:

$$\lim_{w\to\infty}\left\{\sum_{n=1}^{\infty}\left(\frac{1}{\left(n+(-1)+\frac{w-2}{w}\right)\left(n+(-1)+\frac{w-1}{w}\right)\left(n+\frac{1}{w}\right)\left(n+\frac{2}{w}\right)}\right)\right\} = \sum_{n=1}^{\infty}\left(\frac{1}{n^4}\right)$$

$$= \lim_{w\to\infty}\left\{\frac{w^3}{24}\left[\pi\left(4\cot\left(\frac{\pi}{w}\right) - 2\cot\left(\frac{2\pi}{w}\right)\right) - 3w\right]\right\} = \frac{\pi^4}{90}$$

## References

[1] Pietro Mengoli (1650). *Novae quadraturae arithmeticae, seu de additione fractionum*, Bologna.
[2] J.L. Jensen (1915-1916). *An elementary exposition of the theory of the gamma function*, Ann. Math. **17**, 124-166.